\documentclass[12pt,reqno]{amsart}
\usepackage{amsmath, amsthm, amssymb}
\usepackage{mathtools}
\usepackage{amsaddr}
\usepackage[english]{babel}

\usepackage[margin=1in]{geometry}
\usepackage{parskip}
\usepackage[shortlabels]{enumitem}

\usepackage{xcolor}

\usepackage[bookmarks,colorlinks,breaklinks]{hyperref}
\hypersetup{linkcolor=red,citecolor=blue,
  filecolor=dullmagenta,urlcolor=blue}

\usepackage[numbers]{natbib}
\usepackage{cleveref}

\usepackage{graphicx}
\usepackage{tikz}
\usetikzlibrary{matrix, arrows}
\usetikzlibrary{cd}

\makeatletter
\def\thm@space@setup{%
  \thm@preskip=0.5em\thm@postskip=\thm@preskip%
}
\makeatother
\newtheoremstyle{named}{}{}{\\itshape}{}{\bfseries}{.}{.5em}{\thmnote{#3's }#1}
\theoremstyle{named}

\theoremstyle{plain}
\newtheorem{thm}{Theorem}[section]
\newtheorem{prop}[thm]{Proposition}
\newtheorem{lem}[thm]{Lemma}
\newtheorem{cor}[thm]{Corollary}
\theoremstyle{definition}
\newtheorem{defn}[thm]{Definition}

\theoremstyle{remark}
\newtheorem{rmk}[thm]{Remark}
\crefformat{footnote}{#2\footnotemark[#1]#3}

\newcommand{\mathfont}{\mathbf}

\newcommand{\Z}{\mathfont Z}

\newcommand{\Q}{\mathfont Q}
\newcommand{\QQ}{\mathfont Q}
\newcommand{\Qbar}{\overline{\QQ}}

\newcommand{\F}{\mathfont F}
\newcommand{\FF}{\mathfont F}
\newcommand{\Fbar}{\overline{\FF}}
\newcommand{\Fpbar}{\Fbar_p}

\newcommand{\fG}{\mathfrak{G}}
\newcommand{\fH}{\mathfrak{H}}

\newcommand{\fg}{\mathfrak{g}}

\newcommand{\ft}{\mathfrak{t}}

\newcommand{\cC}{\mathcal{C}}

\newcommand{\cO}{\mathcal{O}}

\newcommand{\cT}{\mathcal{T}}
\newcommand{\cU}{\mathcal{U}}

\DeclareFontFamily{OT1}{rsfs}{}
\DeclareFontShape{OT1}{rsfs}{n}{it}{<-> rsfs10}{}
\DeclareMathAlphabet{\mathscr}{OT1}{rsfs}{n}{it}

\newcommand{\into}{\hookrightarrow}

\DeclareMathOperator{\Hom}{Hom}
\DeclareMathOperator{\Ker}{ker}

\DeclareMathOperator{\Gal}{Gal}

\newcommand{\kbar}{\bar{k}}

\DeclareMathOperator{\der}{der}

\newcommand{\Zp}{\mathfont{Z}_p}

\newcommand{\Qp}{\mathfont{Q}_{p}}
\newcommand{\Qpbar}{\overline{\mathfont{Q}}_p}

\newcommand{\rhobar}{{\overline{\rho}}}

\DeclareMathOperator{\Lie}{Lie}

\DeclareMathOperator{\GL}{GL}
\DeclareMathOperator{\SL}{SL}

\newcommand{\Gm}{\mathfont{G}_m}

\DeclareMathOperator{\Image}{Im}

\newcommand{\nlz}[2]{N_{#1}(#2)}

\newcommand{\Dt}{\Delta}

\newcommand{\galabs}[1]{\Gamma_{#1}}
\newcommand{\galres}[2]{\Gamma_{{#1},{#2}}}

\makeatother
\begin{document}

%\subjclass[2010]{11F80}

\title{A note on Galois representations valued in reductive groups with open image \\
}

\author{Shiang Tang}

%\address{150 N University Street \\West Lafayette, IN 47907 \\ United States}

\email{shiangtang1989@gmail.com}

%\keywords{Galois representations, spin groups, inverse Galois problem, automorphic representations}

\maketitle

\begin{abstract}
Let $G$ be a split reductive group with $\dim Z(G) \leq 1$. We show that for any prime $p$ that is large enough relative to $G$, there is a finitely ramified Galois representation $\rho \colon \galabs{\Q} \to G(\Zp)$ with open image. We also show that for any given integer $e$, if the index of irregularity of $p$ is at most $e$ and if $p$ is large enough relative to $G$ and $e$, then there is a Galois representation $\rho \colon \galabs{\Q} \to G(\Zp)$ ramified only at $p$ with open image, generalizing a theorem of Ray \cite{ray:oneprime}. The first type of Galois representation is constructed by lifting a suitable Galois representation into $G(\F_p)$ using a lifting theorem of Fakhruddin--Khare--Patrikis \cite{fkp:reldef}, and the second type of Galois representation is constructed using a variant of the argument in Ray's work \cite{ray:oneprime}. 
\end{abstract}

\section{Introduction}

%\textcolor{red}{compare the first theorem with my ANT paper, say something about geometric gal reps, why do we find such results interesting (galQ and Zp pts of G both profinite, natural to ask blah, galQ has very rich structure), say something about the proofs, FKP does not control ram}

Galois representations arise naturally in algebraic number theory, for example, the $p$-adic Tate module of a rational elliptic curve $E$ gives rise to a continuous representation $\rho \colon \galabs{\Q} \to \GL_2(\Qp)$, where $\galabs{\Q} \coloneqq \Gal(\Qbar/\QQ)$. Moreover, by a result of Serre, $\rho$ has open-image when $E$ is non-CM. Such a Galois representation is part of the rich theory of \textit{geometric Galois representations} (in the sense of Fontaine--Mazur). On the other hand, one does not expect Galois representations $\rho \colon \galabs{\Q} \to \SL_2(\Qp)$ with open image that comes from the cohomology of algebraic varieties or automorphic forms (see, for example \cite[Example 1.4]{tang18}). Thus, for a given reductive algebraic group $G$, it is natural to ask if there is a continuous geometric Galois representation $\rho \colon \galabs{\Q} \to G(\Qp)$ with open image. Several people have constructed examples of this kind (notably for exceptional algebraic groups), for example, \cite{patrikis} and \cite{yun}. On the other hand, for groups like $\SL_2$, it appears to be extremely difficult to disprove the existence of geometric Galois representations with open image. However, if one allows non-geometric Galois representations, a uniform answer can be obtained:  

\begin{thm} [Theorem \ref{thm:galrep}]
Let $G$ be a split reductive group with $\dim Z(G) \leq 1$.
Assume that $p$ is large enough relative to $G$. \footnote{This constant can be made effective, see Remark \ref{rmk:effective bound}.}
Then there is a finitely ramified continuous representation $\rho \colon \galabs{\Q} \to G(\Zp)$ with open image. 
\end{thm}

A similar but weaker result has been obtained in \cite{tang18}, where the author proves the existence of Galois representations $\rho \colon \galabs{\Q} \to G(\Qpbar)$ with \textit{Zariski-dense image}. Note that for a $p$-adic field $E$, a compact, Zariski-dense subgroup of $G(E)$ needs not to be open, unless $G$ is semisimple and $E=\Qp$. On the other hand, working with $\Qp$ (instead of its algebraic closure) imposes a condition on the center of $G$ coming from the structure of $\galabs{\Q}$, see Proposition \ref{prop:dim of center}. We prove the above theorem by lifting a suitable mod $p$ Galois representation $\rhobar \colon \galabs{\Q} \to G(\F_p)$ that has no local obstructions using a lifting theorem of
Fakhruddin--Khare--Patrikis \cite[Theorem 6.21]{fkp:reldef} (and the lifts produced by loc. cit. automatically have open image). Note that we have no control over the ramification loci of the lift $\rho \colon \galabs{\Q} \to G(\Zp)$ due to the nature of the Ramakrishna style lifting argument in loc. cit. In contrast, the next theorem produces open-image Galois representations that ramify only at $p$, assuming a more restrictive condition on $p$. Let $\galres{\Q}{\{p\}}$ be the Galois group of the maximal extension of $\Q$ that is unramified away from $p$.

\begin{thm} [Theorem \ref{thm:oneprime}]
Let $G$ be a split reductive group with $\dim Z(G) \leq 1$. Let $e \geq 0$ and let $p$ be a prime number that is large enough relative to $G$ and $e$ whose index of irregularity is at most $e$.
\footnote{See Theorem \ref{thm:oneprime} for the precise statement on $p$.}
%\begin{enumerate}
    %\item $p$ is larger than $n_{\alpha,\beta}$ for all $\alpha, \beta \in \Phi$, where $n_{\alpha, \beta}$ is a positive integer depending only on the root system $\Phi$ (Proposition \ref{prop:chevalleybasis}).
    %\item $p>1+2N_{e+1}$, where $\{N_k\}$ is a sequence of integers depending only on the root system $\Phi$ (Definition \ref{def:the integers N}). 
    %\item The index of irregularity $e_p$ (Definition \ref{def:index}) is at most $e$. 
%\end{enumerate}
Then there is a continuous representation $\rho \colon \galres{\Q}{\{p\}} \to G(\Zp)$ with open image. 
\end{thm}

This generalizes a theorem of Ray \cite{ray:oneprime}, which studies the $\GL_n$ case. 
Following the argument in loc. cit., we lift a mod $p$ Galois representation valued in a maximal torus of $G$ constructed from the mod $p$ cyclotomic character with \textit{no global obstructions}, which demands certain conditions on the $p$-part of the class group of $\Q(\mu_p)$. On the other hand, to ensure that the image of the lift is open, one modifies any $\Z/p^2$-lift $\rho_2 \colon \galabs{\Q} \to G(\Z/p^2)$ by an appropriate cocycle so that its image satisfies a group-theoretic condition expressed in terms of the Lie algebra $\fg$ (Lemma \ref{lem:root height}), which guarantees that any $\Zp$-lift of $\rho_2$ has open image.  
%Items (2) and (3) in the above theorem come in   

A similar result of Cornut and Ray \cite{cr:iwahori} constructs continuous representations $\rho \colon \galres{\Q}{\{p\}} \to G(\Zp)$ with open image for simple adjoint groups $G$ and \textit{regular primes} $p$ using a completely different method, generalizing work of Greenberg \cite{greenberg}. On the other hand, Maire \cite{maire} constructs continuous representations $\rho \colon \galres{\Q}{\{2,p\}} \to \GL_n(\Zp)$ with open image for every prime $p \geq 3$ (with no regularity condition imposed). 

\begin{rmk}
If the numerators of the
Bernoulli numbers are uniformly random modulo odd primes, then the natural density of primes with index of irregularity $r$ should be $e^{-1/2} \frac{1}{2^r r!}$ (see, for example, \cite{bcs}). In particular, the density of primes with index of irregularity at most $r$ should be at least $1-e^{-1/2}\frac{1}{2^r}$, which approaches 1 rapidly as $r$ increases. 
%\textcolor{red}{say something about the conjectural proportion of $p$ satisfying these conditions with references...}
\end{rmk}

\begin{rmk}
Observe that the above theorems imply the existence of open image Galois representations of $\Gamma_F$ for any number field $F$. In fact, if $\rho \colon \galabs{\Q} \to G(\Zp)$ has open image, then $\rho(\galabs{F})$, being a closed subgroup of finite index, is open. 
\end{rmk}

\textbf{Acknowledgements}. We would like to thank Yves Cornulier, Sean Cotner, and Christian Maire for helpful conversations and comments, and thank the anonymous referee for a careful reading and many comments that have helped us clarify the exposition.

\subsection{Notation}
Let $G$ be a split connected reductive group with derived subgroup $G^{der}$. Let $\fg$ (resp. $\fg^{der}$) be a Lie algebra of $G$ (resp. $G^{der}$). When there is no chance of confusion, we will abuse notation and write $\fg$ (resp. $\fg^{der}$) for $\fg \otimes_{\Z} \F_p$ (resp. $\fg^{der} \otimes_{\Z} \F_p$). 

Let $F$ be a number field. We write $\chi$ for the $p$-adic cyclotomic character and $\overline{\chi}$ for its mod $p$ reduction. Let $\galabs{F} \coloneqq \Gal(\overline{F}/F)$ denote the absolute Galois group of $F$. For any finite set of primes $S$ of $F$, let $\galres{F}{S}$ denote $\Gal(F(S)/F)$, where $F(S)$ is the maximal extension of $F$ inside $\overline{F}$ that is unramified outside the primes in $S$. 

Given a homomophism $\rho \colon \Gamma \to H$ some groups $\Gamma$ and $H$, and an $H$-module $V$, we will write $\rho(V)$ for the associated $\Gamma$-module (we will apply this with $V$ the adjoint representation of an algebraic group).  

\section{Finitely ramified Galois representations with open image}

\subsection{Coxeter elements}
In this section, we review the notion of Coxeter elements, for more details, see \cite[\S 10.1]{bhkt}. 
Let $G$ be a split \textit{simple, simply-connected} group with root system $\Phi=\Phi(G,T)$. Let $\Dt=\{\alpha_1, \cdots, \alpha_r\} \subset \Phi$ be a set of simple roots. We write $W=W(G,T)=N_G(T)/T$ for the Weyl group of $G$. We call an element $w \in W$ a \textit{Coxeter element} if it is conjugate in $W$ to an element of the form $w_{\alpha_1} \cdots w_{\alpha_r}$, where $w_{\alpha_i}$ is the simple reflection corresponding to $\alpha_i$. There is a unique conjugacy class of Coxeter elements in $W$. Their common order $h$ is called the \textit{Coxeter number} of $G$. 

\begin{lem} \label{lem:lifting W}
Let $w$ be a Coxeter element. Then there is an element $\tilde w \in \nlz{G}{T}(\Z)$ lifting $w$. Its order $\tilde h$ depends only on $G$.
%There exists an element $\tilde w_i \in \nlz{G_i(\Z)}{T_i}$ whose image in $W(G_i,T_i)$ equals $w_i$ such that the order of $w_i$ is either $h_i$ or $2h_i$. Here $h_i$ denotes the Coxeter number of $G_i$. 
\end{lem}
\begin{proof}
By \cite{tits}, there is a finite subgroup $\cT \subset \nlz{G}{T}(\Z)$ which is isomorphic to the extension of $W(G,T)$ by a subgroup of $T(\Z)$. In particular, any element in $W$ lifts to a finite order element in $\nlz{G}{T}(\Z)$.
That the order of $\tilde w$ depends only on $G$ follows from \cite[Proposition 10.2, (iii)]{bhkt}. 
\end{proof}

\begin{prop} \label{prop:coxeter irrd}
Let $k$ be a field of characteristic $p$. Assume that $p>2h-2$. Let $\Gamma \subset G(k)$ be a subgroup. Assume that
\begin{enumerate}
    \item $\Gamma \subset N_G(T)(k)$.
    \item The image of $\Gamma$ in $W(G,T)$ is the cyclic group generated by a Coxeter element.
\end{enumerate}
Then $\Gamma$ is $G$-irreducible, i.e. $\Gamma$ is not contained in any proper parabolic subgroup of $G(\kbar)$. 
\end{prop}
\begin{proof}
This follows from the proof of \cite[Proposition 10.7, (i)]{bhkt}.
\end{proof}

\subsection{Lifting Galois representations}
In this section, we assume that $G=G_1 \times \cdots \times G_n$ is a direct product of simple, simply connected groups. For $1 \leq i \leq n$, let $T_i$ be a maximal torus of $G_i$ and let $T=T_1 \times \cdots \times T_n$. For each $i$, let $\tilde w_i \in \nlz{G_i}{T_i}(\Z)$ be as in Lemma \ref{lem:lifting W}. Write $\tilde h_i$ for its order. % and let $\tilde w_i \in \nlz{G_i(\Z)}{T_i}$ be a lift of $w_i$. 

\begin{prop} \label{prop:igp}
Let $\fG \subset N_G(T)(\F_p)$ be the group generated by $T(\F_p)$ and the elements $\tilde w_i$ for $1 \leq i \leq n$. Assume that $p>c(\tilde h_1, \cdots, \tilde h_n)$ (a constant depending only on $\tilde h_i$ for $1 \leq i \leq n$). Then there is a finite Galois extension $M/\Q$ whose Galois group is isomorphic to $\fG$ in which $p$ is unramified. 
\end{prop}
\begin{proof}
This is an easy consequence of \cite[Theorem 2.5]{tang18}. Let $\fH$ be the group generated by $\tilde w_i$ for $1 \leq i \leq n$, then $\fG$ is a quotient of the semidirect product of $\fH$ and $T(\F_p)$. %Let $K/\Q$ be a Galois extension whose Galois group is isomorphic to $\fH \cong \prod_i \Z/\tilde h_i$: 
For each $i$, let $p_i$ be an odd prime for which $p_i \equiv 1 \pmod{\tilde h_i}$, and let $K_i$ be the fixed field of the unique subgroup of $\Gal(\Q(\mu_{p_i})/\Q)$ of index $\tilde h_i$, and let $K=K_1 \cdots K_n$. Then $\Gal(K/\Q) \cong \prod_i \Z/\tilde h_i \cong \fH$ as long as the primes $p_i$ are chosen to be distinct. Let $c(\tilde h_1, \cdots, \tilde h_n)$ be the smallest possible value of $\max \{p_1, \cdots, p_n\}$.  
Then loc. cit. implies that if $p$ is unramified in $K$ (in particular, if $p>c(\tilde h_1, \cdots, \tilde h_n)$), then there is a number field $M$ with claimed properties.  
\end{proof}

\begin{rmk} \label{rmk:primes}
It is possible to work out an explicit bound for $p$ in Proposition \ref{prop:igp}. In fact, \cite{tv:leastprime} shows that the least prime congruent to 1 modulo $n$ is at most $2^{\phi(n)+1}-1$. 
\end{rmk}

Let $\rhobar \colon \galabs{\Q} \to G(\F_p)$ be the mod $p$ Galois representation associated to the extension $M/\Q$ constructed above. 

\begin{prop} \label{prop:lift}
Assume that $p$ is greater than $\max \{ 2h_i-2 | 1 \leq i \leq n \}$, $\max \{ \tilde h_i | 1 \leq i \leq n \}$, $c(\tilde h_1, \cdots, \tilde h_n)$, and $c_G$, where $c_G$ is a constant depending only on $G$ in \cite[Theorem 6.11]{fkp:reldef}.  
%(\textcolor{red}{coming from either the previous props and FKP's Thm 6.21. The constant in FKP seems hard to extract but in our special case this might be possible }). 
Then $\rhobar$ lifts to a finitely ramified continuous representation
$\rho \colon \galabs{\Q} \to G(\Zp)$ whose image contains $\widehat G(\Zp)$. 
\end{prop}
\begin{proof}
We apply \cite[Theorem 6.21]{fkp:reldef}. We will check its assumptions and explain why the lift can be chosen to have $\Zp$-coefficients. First note that the field $\widetilde{F}$ in the statement of loc. cit. equals $\QQ$ (since $G$ is connected) and $[\QQ(\mu_p):\QQ]=p-1$.
The first item of loc. cit. holds since the projection of $\fG$ to $G_i(\Fpbar)$ is $G_i$-irreducible by Proposition \ref{prop:coxeter irrd}.
We claim that for all finite primes $v$, there exists a formally smooth deformation condition for $\rhobar|_{\galabs{\Q_v}}$. In fact, for $v=p$, $H^2(\galabs{\Q_p}, \rhobar(\fg))=H^0(\galabs{\Q_p}, \rhobar(\fg)(1)) \subset H^0(I_{\Qp}, \rhobar(\fg)(1))=0$, where the first equality follows from local duality, and the last equality follows from the fact that $\rhobar(I_{\Qp})=1$ (which holds since $p$ is unramified in the fixed field of $\rhobar$). For $v \neq p$, since $\rhobar(I_{\Q_v}) \subset \fG$ is prime to $p$ (by the construction of $\fG$), we can take the (formally smooth) minimal prime to $p$ deformation condition for $\rhobar|_{\galabs{\Q_v}}$ (see \cite[\S 4.4]{patrikis}). Thus, the second item of \cite[Theorem 6.21]{fkp:reldef} holds. By loc. cit., $\rhobar$ lifts to a finitely ramified representation $\rho \colon \galabs{\Q} \to G(\cO)$ whose image contains $\widehat{G}(\cO)$ for a finite extension $\cO$ of $\Zp$. Lastly, note that we may take $\cO=\Zp$ in our case since a formally smooth deformation condition exists at every finite place, see the third item of \cite[Remark 1.3]{fkp:reldef}. 
\end{proof}

\subsection{The general case}

%\begin{lem} \label{lem:p adic lie group}
%Let $G$ be a semisimple algebraic group, then any compact, Zariski-dense subgroup of $G(\Qp)$ is open.
%\end{lem}

%\begin{proof}
%Let $C$ be such a group, which is $\Qp$-analytic. The Lie algebra $\Lie C$ is $C$-invariant for the adjoint representation, hence is $G$-invariant by Zariski density, so $\Lie C$ is an ideal of the semisimple Lie algebra $\Lie G$. So we can write $G=HK$ for commuting semisimple algebraic subgroups $H$ and $K$ such that $H \cap K$ is finite and $\Lie H = \Lie C$. Then $C$ has a finite index subgroup contained in $H(\Qp)$, and hence the image of $C$ in $G/H$ is finite, which implies $G=H$ by Zariski density. So $\Lie C=\Lie G$ and thus $C$ is open.  
%\end{proof}

\begin{lem} \label{lem:submersion}
Let $\widetilde G$ and $G$ be algebraic groups defined over $\Qp$ and let $\widetilde{G} \to G$ be an isogeny. Then the induced map $\widetilde G(\Qp) \to G(\Qp)$ is open.
\end{lem}

\begin{proof}
First note that if $f \colon X \to Y$ is a submersion of (real or $p$-adic) manifolds, then $f$ is open by the local structure theorem for submersions (see for example, \cite[Part II, Ch. III, Section 10]{serre:lie}). Moreover, the algebraic and analytic differentials are compatible, so if $X,Y$ are varieties over $\Qp$ with $\dim X \geq \dim Y$ and if $f \colon X \to Y$ is a smooth morphism, then the induced map $X(\Qp) \to Y(\Qp)$ is a submersion, and hence open. Now if $\widetilde G \to G$ is an isogeny of algebraic groups defined over $\Qp$, then it is smooth, so $\widetilde G(\Qp) \to G(\Qp)$ is open by the above. 
%First a general thing: if f: X --> Y is a submersion of (real or p-adic) manifolds, then f is open, as follows from the local structure theorem for submersions. Moreover, the algebraic and analytic differentials are compatible, so if f: X --> Y is a smooth morphism of varieties over Q_p then the induced map X(Q_p) --> Y(Q_p) is a submersion, hence open. If G' --> G is an isogeny of Q_p-group schemes, then it is smooth, so G'(Q_p) --> G(Q_p) is open. As G'(Z_p) is open in G'(Q_p), its image in G(Q_p), and hence also G(Z_p), is open. 
\end{proof}

We thank Sean Cotner for pointing out the above lemma.

\begin{thm} \label{thm:galrep}
Let $G$ be a split reductive group with $\dim Z(G) \leq 1$.
Assume that $p$ is large enough relative to $G$. Then there is a finitely ramified continuous representation $\rho \colon \galabs{\Q} \to G(\Zp)$ with open image. 
\end{thm}
\begin{proof}
First assume that $G$ is semisimple, so there are simple, simply-connected groups $G_1, \cdots, G_n$, equipped with an isogeny $\widetilde{G}=G_1 \times \cdots \times G_n \to G$. Construct $\rhobar \colon \galabs{\Q} \to \widetilde{G}(\F_p)$ as in the previous section and apply Proposition \ref{prop:lift}, we obtain a finitely ramified representation $\widetilde\rho \colon \galabs{\Q} \to \widetilde G(\Zp)$ with open image. Projecting down to $G$, we obtain a Galois representation with open image by Lemma \ref{lem:submersion}. If $G$ is reductive with $\dim Z(G)=1$, there is a canonical isogeny $G^{der} \times Z(G)^0 \to G$. The above gives an open image Galois representation into $G^{der}(\Zp)$. On the other hand, the cyclotomic character $\chi \colon \galres{\Q}{p } \to \Gm(\Zp) = Z(G)^0(\Zp)$ has open image. Thus by Lemma \ref{lem:submersion}, their product gives a Galois representation into $G(\Zp)$ with open image.
%\textcolor{red}{may need to assume that $p$ does not divide the kernel of these isogenies}
\end{proof}

\begin{rmk} \label{rmk:effective bound}
The lower bound for $p$ in Theorem \ref{thm:galrep} can be made effective: by its proof, this bound is the maximum of the four constants in Proposition \ref{prop:lift} associated to the simply-connected cover of $G$. By Remark \ref{rmk:primes}, the constant $c(\tilde h_1, \cdots, \tilde h_n)$ can be bounded by an explicit formula, and by \cite[Remark 6.17]{fkp:reldef}, the constant $c_G$ can be made effective as well. 
\end{rmk}

\begin{prop} \label{prop:dim of center}
Let $G$ be a split reductive group. Suppose that there exists a continuous representation $\rho \colon \galabs{\Q} \to G(\Zp)$ with open image. Then $\dim Z(G) \leq 1$.
\end{prop}
\begin{proof}
Let $S(G)=G/G^{der}$. It is a split torus with $\dim Z(G) = \dim S(G) =: r$. Since $\rho$ has open image, so does $\rho \pmod{G^{der}} \colon \galabs{\Q} \to S(G)(\Zp)=\Gm(\Zp)^r$.
In fact, $\Image \rho \cap Z(G)^0(\Qp)$ is an open subgroup of $Z(G)^0(\Qp)$, which maps to an open subgroup of $S(G)^0(\Qp)$ under the canonical isogeny of tori $Z(G)^0 \to S(G)^0$ by Lemma \ref{lem:submersion}. 
Since $\Q$ has a unique $\Zp$-extension, this forces $r$ to be at most 1. 
\end{proof}

\section{Galois representations ramified at one prime with open image}
Suppose that $G$ is a split reductive group.
Given a continuous representation $\rhobar \colon \galres{\Q}{\{p\}} \to G(\F_p)$. %Set $\overline{\mu} \coloneqq \rhobar \pmod{G^{der}} \colon \galres{\Q}{\{p\}} \to G/G^{der}(\FF_p)$ and fix a lift $\mu \colon \galres{\Q}{\{p\}} \to G/G^{der}(\Zp)$ which is unramified away from $p$ and (\textcolor{red}{something else}). For the rest of this paper, we require all lifts $\rho \colon \galres{\Q}{\{p\}} \to G(R)$ of $\rhobar$ (with $R$ a $\Zp$-algebra) to satisfy $\rho \pmod{G^{der}} = \mu$. 
Suppose that $\rho \colon \galres{\Q}{\{p\}} \to G(\Zp)$ is a continuous lift of $\rhobar$. For $m \geq 1$, set $\rho_m$ to be the mod-$p^m$ reduction of $\rho$ mod $p^m$. 

The following fact is standard, see \cite[\S 3.5]{tilouine}:

\begin{lem} \label{lem:expmap}
There is a natural group isomorphism
\[ \exp \colon \fg \otimes_{\F_p} p^m\Z/p^{m+1}\Z 
\xrightarrow{\sim} \Ker (G(\Z/p^{m+1}) \to G(\Z/p^m))
\]
%and the same holds with $G$ (resp. $\fg$) replaced by $G^{der}$ (resp. $\fg^{der}$).
\end{lem}

\begin{defn}
For $m \geq 1$, set $\Phi_m(\rho) \coloneqq \rho_{m+1}(\ker \rho_m) \subset G(\Z/p^{m+1})$.
\end{defn}

The following lemma follows immediately from the above, we omit the proof.

\begin{lem} \label{lem:phimrho}
$\Phi_m(\rho)$ may be identified as a submodule of $\rhobar(\fg)$: for $g \in \ker \rho_m$, $\rho_{m+1}(g)=\exp(p^m v)$ for a unique $v \in \rhobar(\fg)$, and we identify $\rho_{m+1}(g)$ with this $v$.
\end{lem}

\begin{lem}\label{lem:lie bracket}
With the identification in Lemma \ref{lem:phimrho}, for $l,m \geq 1$, $[\Phi_l(\rho),\Phi_m(\rho)] \subset \Phi_{l+m}(\rho)$, where $[,]$ is the Lie bracket of $\fg$.
\end{lem}
\begin{proof}
%By \cite[Lemma 2.8]{ray:oneprime}, the lemma holds for $G=\GL_n$. 
Fix a faithful representation $i \colon G \into \GL_n$ defined over $\Z$ for some integer $n$, which induces a map on the Lie algebras as well. Note that $i\Phi_m(\rho)=\Phi_m(i\rho)$. It follows that $i[\Phi_l(\rho),\Phi_m(\rho)]=[i\Phi_l(\rho),i\Phi_m(\rho)]=[\Phi_l(i\rho),\Phi_m(i\rho)] \subset \Phi_{l+m}(i\rho)=i\Phi_{l+m}(\rho)$ where the $\subset$ in the middle follows from \cite[Lemma 2.8]{ray:oneprime}. So $[\Phi_l(\rho),\Phi_m(\rho)] \subset \Phi_{l+m}(\rho)$.
\end{proof}

\begin{lem}
Let $\rho \colon \galres{\Q}{\{p\}} \to G(\Zp)$ be a continuous representation lifting $\rhobar$. Let $m \geq 1$ be such that $\Phi_m(\rho)$ contains $\rhobar(\fg^{der})$. Then we have
\begin{enumerate}
    \item $\Phi_k(\rho)$ contains $\rhobar(\fg^{der})$ for $k $ divisible by $m$. 
    \item The image of $\rho$ contains $\cU_m \coloneqq \ker G^{der}(\Zp) \to G^{der}(\Zp/p^m)$.
\end{enumerate}
\end{lem}
\begin{proof}
(1) follows from the identity $[\rhobar(\fg^{der}),\rhobar(\fg^{der})]=\rhobar(\fg^{der})$ and Lemma \ref{lem:lie bracket}. Let $H$ be the image of $\rho$. %Since $\rho$ is continuous and $\galres{\Q}{\{p\}} $ is compact, it follows that $H$ is closed. Then we may identify $H$ with the inverse limit $\inverselimit_k H_k$ where $H_k$ is the projection of $H$ to $G(\Z/p^k)$. 
(1) implies that for infinitely many $k$, $H \pmod {p^k}$ contains $\cU_m \pmod{p^k}$. %and hence $H=\inverselimit_k H_k$ contains $\cU_m$. 
%Furthermore, if $k \geq N_1$, then (taking $P=H \cap \cU_m$) Lemma \ref{fkp lemma} implies that $H $ contains $\cU_m$. 
If moreover, $k \geq N_1$, then we can apply the lemma below (with $P=H \cap \cU_m$) to conclude that that $H $ contains $\cU_m$. 
\end{proof}

\begin{lem}\label{fkp lemma}
Let $m$ be a positive integer.
There is an integer $N_1 \geq m$ depending only on $m$ and $G^{\der}$ such that if $k \geq N_1$,
any closed subgroup $P$ of $\cU_m$ whose reduction modulo $p^k$ equals $\ker G^{\der}(\Zp/p^k) \to G^{\der}(\Zp/p^m)$ must in fact equal $\cU_m$. 
\end{lem}

\begin{proof}
This is essentially \cite[Lemma 6.15]{fkp:reldef}, which proves the case when $m=1$. The argument trivially generalizes to arbitrary $m$.   
\end{proof}

For the rest of this section, we fix a maximal split torus $T$ of $G$. Let $\Phi=\Phi(G,T)$ be the associated root system. Fix a set of simple roots $\Dt \subset \Phi$. Let $\Phi^+$ denote the corresponding set of positive roots. For any $\alpha \in \Phi$, write $\alpha=\sum_{\beta \in \Dt} n_{\beta} \beta$ and define the height function
\[
\mathrm{ht}(\alpha) \coloneqq \sum_{\beta \in \Dt} n_{\beta}. 
\]
The following proposition is well-known, see \cite[Chapter 1]{steinberg}.
\begin{prop} \label{prop:chevalleybasis}
There exist elements $H_{\alpha} \in \ft=\Lie T$ and $X_{\alpha} \in \fg^{der}=\Lie G^{der}$ for $\alpha \in \Phi$, such that 
the elements $H_{\alpha}$ for $\alpha \in \Dt$ and $X_{\alpha}$ for $\alpha \in \Phi$
%$\{ H_{\alpha}, \alpha \in \Dt, X_{\alpha}, \alpha \in \Phi \}$ 
form an integral basis for $\fg$ satisfying the relations below: 
\begin{itemize}
    \item $[H_{\alpha}, H_{\beta}]=0$.
    \item $[H_{\beta}, X_{\alpha}]=\alpha(H_{\beta}) X_{\alpha}$ with $\alpha(H_{\beta}) \in \Z$. %is an integer depending on $\alpha$ and $\beta$ (called a Cartan integer).
    \item $[X_{\alpha}, X_{-\alpha}]=H_{\alpha}$, and $H_{\alpha}$ is an integral combination of the $H_{\beta}$ with $\beta \in \Dt$.
    \item $[X_{\alpha}, X_{\beta}]=N_{\alpha,\beta} X_{\alpha+\beta}$ with $N_{\alpha,\beta} \in \Z-\{0\}$ if $\alpha+\beta \in \Phi$.
    \item $[X_{\alpha}, X_{\beta}]=0$ if $\alpha+\beta \neq 0$ and $\alpha+\beta \notin \Phi$.
\end{itemize}
Such a basis is called a \textbf{Chevalley basis}. It is unique up to sign changes and automorphisms. 
Moreover, the constant $n_{\alpha,\beta} \coloneqq |N_{\alpha,\beta}|$ can be described purely in terms of the root system $\Phi$.  
\end{prop}

\begin{lem} \label{lem:root height}
Let $\rho \colon \galres{\Q}{\{p\}} \to G(\Zp)$ be a continuous representation lifting $\rhobar$.
Assume that $p$ is larger than $n_{\alpha,\beta}$ for all $\alpha, \beta \in \Phi$.
Assume that $\Phi_1(\rho)$ contains an element $H \in \ft \subset \rhobar(\fg)$ such that $\alpha(H)$ is a nonzero element in $\F_p$ for all $\alpha \in \Phi$. Furthermore, assume that it contains $X_{\alpha}$ for all $\alpha \in \Phi$ with $\mathrm{ht}(\alpha)$ odd. Then we have
\begin{enumerate}
    \item $\Phi_{4}(\rho)$ contains $\rhobar(\fg^{der})$.
    \item The image of $\rho$ contains $\cU_{4}$.
\end{enumerate}
\end{lem}
\begin{proof}
In the course of the proof, we will use Lemma \ref{lem:lie bracket} several times without reference.
Let $\alpha \in \Phi^+$ be a root with $\mathrm{ht}(\alpha)=n$ even. Write $\alpha=\beta+\gamma$ for some $\beta \in \Dt$ and $\gamma \in \Phi^+$. Then $\mathrm{ht}(\beta)=1$ and $\mathrm{ht}(\gamma)=n-1$, both are odd. The relation $N_{\beta,\gamma} X_{\beta+\gamma} = [X_{\beta}, X_{\gamma}]$ and the assumption on $p$ imply that $X_{\alpha}=X_{\beta+\gamma} \in \Phi_2(\rho)$. A similar argument gives that $X_{\alpha} \in \Phi_2(\rho)$ for $\alpha \in \Phi^-$ with even height.  
Now let $\alpha \in \Phi$ be a root with $\mathrm{ht}(\alpha)$ odd. The relation $[H, X_{\alpha}]=\alpha(H) X_{\alpha}$ and the assumption on $H$ imply that $X_{\alpha} \in \Phi_2(\rho)$. 

Since $[\Phi_1(\rho),\Phi_2(\rho)] \subset \Phi_3(\rho)$, the relation $[H, X_{\alpha}]=\alpha(H) X_{\alpha}$ implies that $X_{\alpha} \in \Phi_3(\rho)$ for all $\alpha \in \Phi$. One more iteration of the same kind implies that $X_{\alpha} \in \Phi_4(\rho)$ for all $\alpha \in \Phi$.
That $[\Phi_2(\rho),\Phi_2(\rho)] \subset \Phi_4(\rho)$ and the relation $[X_{\alpha}, X_{-\alpha}]=H_{\alpha}$ imply that $H_{\alpha} \in \Phi_4(\rho)$ for all $\alpha \in \Phi$. Thus, $\Phi_4(\rho)$ contains $\rhobar(\fg^{der})$. The previous lemma now implies our claims.
\end{proof}

Let $\lambda \in X_*(T)$ be a cocharacter. We will impose conditions on $\lambda$ later.
Let 
\[ \langle , \rangle \colon X^*(T) \times X_*(T) \to \Hom(\Gm, \Gm) = \Z \]
denote the canonical pairing. 
Let $\rhobar \colon \galres{\QQ}{\{p\}} \to T(\FF_p)$ be $\lambda \circ \overline{\chi}$ (recall that $\overline{\chi}$ is the mod $p$ cyclotomic character). 

%\textcolor{red}{define $\cC(?)$ and index of irregularity somewhere}.

Let $\cC$ denote the mod $p$ class group of $\Q(\mu_p)$, i.e. $\cC=\mathrm{Cl}(\Q(\mu_p)) \otimes \F_p$. It has a natural action of $\Gal(\Q(\mu_p)/\Q)$ and $\cC$ decomposes into eigenspaces 
\[
\cC=\bigoplus_{0 \leq i \leq p-2} \cC(\overline{\chi}^i)
\]
where $\cC(\overline{\chi}^i)=\{ x \in \cC | g \cdot x = \overline{\chi}^i(g) x \}$. 

\begin{defn} \label{def:index}
The index of irregularity $e_p$ of a prime $p$ is the number of eigenspaces $\cC(\overline{\chi}^i)$ which are nonzero. If $e_p=0$, we say $p$ is regular. 
\end{defn}

\begin{thm} \label{thm:key theorem}
Let $\lambda$ and $\rhobar$ be as above. Assume that
\begin{enumerate}
    \item $p$ is larger than $n_{\alpha,\beta}$ for all $\alpha, \beta \in \Phi$.
    \item $0<|\langle \alpha, \lambda \rangle|<p-1$ for all $\alpha \in \Phi$.
    \item $\langle \alpha, \lambda \rangle$ is odd for all $\alpha \in \Dt$.
    \item The characters $\overline{\chi}^{\langle \alpha, \lambda \rangle} $ for all $\alpha \in \Phi$ are all distinct and are not equal to $\overline{\chi}$.  
    \item $\cC(\overline{\chi}^{p-\langle \alpha, \lambda \rangle})=0$ for all $\alpha \in \Phi$.
\end{enumerate}
Then $\rhobar$ admits a continuous lift $\rho \colon \galres{\Q}{\{p\}} \to G(\Zp)$ whose image contains $\cU_4$. 
\end{thm}

\begin{proof}
We follow the argument of \cite[Theorem 3.3]{ray:oneprime} closely.
First, we have $H^2(\galres{\Q}{\{p\}}, \rhobar(\fg))=0$. This follows from the assumptions that $\overline{\chi}^{\langle \alpha, \lambda \rangle} \neq \overline{\chi}$, $\cC(\overline{\chi}^{p-\langle \alpha, \lambda \rangle})=0$, and the local and global duality theorems. We refer the reader to the first part of the proof of \cite[Theorem 3.3]{ray:oneprime} for details (the argument in loc. cit. is for $\GL_n$ but it trivially generalizes to $G$). 

Let $\chi_2$ be $\chi$ mod $p^2$ and let $\rho_2' \coloneqq \lambda \circ \chi_2$. Let $\alpha \in \Phi$ be a root with odd height. Then by assumption, $\langle \alpha, \lambda \rangle$ is odd, and hence $H^0(\galabs{\Q_{\infty}}, \FF_p(\overline{\chi}^{\langle \alpha, \lambda \rangle}))=0$ and $H^0(\galres{\Q}{\{p\}}, \FF_p(\overline{\chi}^{\langle \alpha, \lambda \rangle}))=0$. By the previous paragraph, $H^2(\galres{\Q}{\{p\}}, \FF_p(\overline{\chi}^{\langle \alpha, \lambda \rangle}))=0$. It follows from the global Euler characteristic formula that $H^1(\galres{\Q}{\{p\}}, \FF_p(\overline{\chi}^{\langle \alpha, \lambda \rangle}))$ is 1-dimensional. Let $f_{\alpha}$ be a generator of $H^1(\galres{\Q}{\{p\}}, \FF_p(\overline{\chi}^{\langle \alpha, \lambda \rangle}))$ and let $F \in H^1(\galres{\Q}{\{p\}}, \rhobar(\fg))$ be the sum of all $f_{\alpha}$ with $\alpha$ ranging over roots in $\Phi$ with odd height. Let 
\[ \rho_2 \coloneqq \exp(pF) \cdot \rho_2'=\exp(pF) \cdot (\lambda \circ \chi_2). \] 
As $H^2(\galres{\Q}{\{p\}}, \rhobar(\fg))=0$, $\rho_2 $ lifts to a characteristic zero representation $\rho \colon \galres{\Q}{\{p\}} \to G(\Zp)$. 

We want to show that the image of $\rho$ contains $\cU_4$. By Lemma \ref{lem:root height}, it suffices to show that $\Phi_1(\rho)$ contains
\begin{itemize}
    \item $X_{\alpha}$ for all $\alpha \in \Phi$ with $\mathrm{ht}(\alpha)$ odd,
    \item an element $H \in \ft$ such that $\alpha(H)$ is a nonzero element in $\F_p$ for all $\alpha \in \Phi$. 
\end{itemize}
Since the image of $\rhobar$ is prime to $p$, any Galois submodule $M$ of $\rhobar(\fg)$ decomposes into 
\[
M=M_1 \oplus \big ( \bigoplus_{\alpha \in \Phi} M_{\overline{\chi}^{\langle \alpha, \lambda \rangle}} \big )
\]
where $M_1$ is the $\galabs{\Q}$-invariant submodule and $M_{\overline{\chi}^{\langle \alpha, \lambda \rangle}}$ is the $\overline{\chi}^{\langle \alpha, \lambda \rangle}$-eigenspace. Since the characters $\overline{\chi}^{\langle \alpha, \lambda \rangle} $ for all $\alpha \in \Phi$ are all distinct and are nontrivial (by assumption), the above decomposition makes sense and $M_{\overline{\chi}^{\langle \alpha, \lambda \rangle}}$, if nonzero, is the 1-dimensional space generated by $X_{\alpha}$. Note that as $H^1(\Gal(\Q(\mu_p)/\Q), \FF_p(\overline{\chi}^{\langle \alpha, \lambda \rangle}))=0$, it follows from the inflation-restriction sequence that for any root $\alpha$ with odd height, the restriction of $f_{\alpha}$ to $\galabs{\Q(\mu_p)}$ is nonzero. Hence, there exists $g \in \ker \rhobar$ such that $f_{\alpha}(g) \neq 0$, and so the element $\rho_2(g) \in \Phi_1(\rho)$ has nonzero $X_{\alpha}$-component. It follows from the above decomposition with $M=\Phi_1(\rho)$ that $X_{\alpha} \in \Phi_1(\rho)$ for all $\alpha \in \Phi$ with odd height. 

Finally, note that the cyclotomic character $\chi$ induces an isomorphism $\chi \colon \Gal(\Q(\mu_{p^\infty})/\Q(\mu_p)) \xrightarrow{\sim} 1+p\Zp$. Let $\gamma \in \galabs{\Q(\mu_p)}$ be chosen 
such that $\chi(\gamma)=1+p$. Then $\rho_2(\gamma)=\exp(pF(\gamma)) \cdot \lambda(\chi_2(\gamma))=\exp(pF(\gamma)) \cdot \lambda(1+p) \in \Phi_1(\rho)$. Let $H \in \ft \subset \rhobar(\fg)$ be the element
such that $\exp(pH)=\lambda(1+p) $. Then $H \in \Phi_1(\rho)$ and for any root $\alpha$, $\alpha(H)=\langle \alpha, \lambda \rangle$ is nonzero mod $p$ by Assumption (2). 
\end{proof}

Let $r$ be the rank of $\Phi$ and let $\Dt=\{ \alpha_1, \cdots, \alpha_r \}$. Let $c_1, \cdots, c_r$ be positive integers such that $\tilde \alpha = \sum c_i \alpha_i$ is the highest root. Let $\lambda_1, \cdots, \lambda_r$ be cocharacters such that $\langle \alpha_i, \lambda_j \rangle =\delta_{ij}$. 

\begin{defn} \label{def:the integers N}
Define an strictly increasing sequence of integers $N_0, N_1, N_2, N_3, \cdots$ as follows: $N_0=1$, 
\[
N_{k+1}=c_1(N_k^*+2)+c_2(N_k^*+4)+\cdots+c_r(N_k^*+2r)
\]
where $N_k^*=N_k$ if $N_k$ is odd, and $N_k^*=N_k+1$ if $N_k$ is even. 
\end{defn}

\begin{cor} \label{cor:galrep one prime}
Let $e \geq 0$ and let $p$ be a prime number such that
\begin{enumerate}
    \item $p$ is larger than $n_{\alpha,\beta}$ for all $\alpha, \beta \in \Phi$.
    \item $p>1+2N_{e+1}$.
    \item The index of irregularity $e_p$ is at most $e$.
\end{enumerate}
Then there is a continuous representation $\rho \colon \galres{\Q}{\{p\}} \to G(\Zp)$ whose image contains $\cU_4$. 
\end{cor}
\begin{proof}
We only need to show that for any prime $p$ satisfying the above conditions, there is a cocharacter $\lambda$ satisfying the assumptions in Theorem \ref{thm:key theorem}. Let $A \subset \Z/(p-1)$ be defined by $n \in A$ if and only if at least one of $\cC(\overline{\chi}^{p+n})$ or $\cC(\overline{\chi}^{p-n})$ is nonzero. Then since $e_p \leq e$, we have $|A| \leq 2e$.
It suffices to show that there is a cocharacter $\lambda$ such that $\lambda $ satisfies Theorem \ref{thm:key theorem}, (2)-(4) and $\langle \alpha, \lambda \rangle \notin A$ for all $\alpha \in \Phi^+$. We may assume that $e_p=e \geq 1$ and write the lease positive representatives of the elements in $A$
in ascending order
\[
0 \leq a_1 < \cdots < a_e \leq p-1-a_e < \cdots < p-1-a_1 \leq p-1.
\]
We prove by induction that if $p>1+2N_{e+1}$, then there exists a cocharacter $\lambda$ such that 
\begin{itemize}
    \item For $1 \leq j \leq r$, let $x_j \coloneqq \langle \alpha_j, \lambda \rangle$. Then $x_j>1$ is odd and $x_1, \cdots, x_r$ is an arithmetic progression with a common difference of 2. 
    \item $S \coloneqq \{\langle \alpha, \lambda \rangle| \alpha \in \Phi^+ \}$ falls in between $a_i$ and $a_{i+1}$ for some $0 \leq i \leq e-1$ (set $a_0=0$) or in between $a_e$ and $p-1-a_e$. 
    \item Moreover, $S$ can be made so that $\max S < a_e$ unless $a_i \leq N_i$ for all $i$ with $1 \leq i \leq e$.
\end{itemize}
Granted this, it clearly implies Theorem \ref{thm:key theorem}, (2)-(5), and so if we further require $p$ to be larger than $n_{\alpha,\beta}$ for all $\alpha, \beta \in \Phi$, then we obtain a desired lift. 

It remains to prove the claim. First suppose that $e=1$ and $p>1+2N_2$. If $a_1>N_1=3c_1+5c_2+\cdots +(2r+1)c_r$, we take $\lambda=3\lambda_1+5\lambda_2+\cdots +(2r+1)\lambda_r$,
then since $\tilde\alpha= c_1\alpha_1+c_2\alpha_2+\cdots +c_r\alpha_r$ is the highest root, for any $\alpha \in \Phi^+$,
$0<\langle \alpha, \lambda \rangle \leq \langle \tilde\alpha, \lambda \rangle=N_1<a_1$. If $a_1 \leq N_1$, we take $\lambda=(N_1^*+2)\lambda_1+\cdots +(N_1^*+2r)\lambda_r$, then for any $\alpha \in \Phi^+$, $a_1 \leq N_1<\langle \alpha, \lambda \rangle \leq \langle \tilde\alpha, \lambda \rangle = N_2 < p-1-N_1 \leq p-1-a_1$, where the second last inequality holds since $p>1+2N_2>1+N_1+N_2$. Thus the claim holds for $e=1$. Assume the claim holds for $e$, and consider a sequence
\[
(0 \leq ) a_1 < \cdots < a_e < a_{e+1} \leq p-1-a_{e+1} < \cdots <p-1-a_1 (\leq p-1).
\]
Let $p$ be a prime with $p>1+2N_{e+2}$. 
If at least one of the $a_i $ with $1 \leq i \leq e$ is greater than $N_i$, the induction hypothesis (applied to the sequence with the $a_{e+1}$ terms removed) implies that there is a cocharacter $\lambda$ satisfying the properties in the claim with $\max S<a_e$ and we are done. (Note that without the condition that $\max S < a_e$, the $S$ provided by induction could intersect with $\{a_{e+1}, p-1-a_{e+1}\}$.) If $a_i \leq N_i$ for all $1 \leq i \leq e$ but $a_{e+1}>N_{e+1}$, we take $\lambda=(N_e^*+2)\lambda_1+\cdots +(N_e^*+2r)\lambda_r$, then for any $\alpha \in \Phi^+$, $a_e \leq N_e<\langle \alpha, \lambda \rangle \leq \langle \tilde\alpha, \lambda \rangle = N_{e+1} < a_{e+1}$ and we are done. If $a_i \leq N_i$ for all $1 \leq i \leq e+1$, we take $\lambda=(N_{e+1}^*+2)\lambda_1+\cdots +(N_{e+1}^*+2r)\lambda_r$, then for any $\alpha \in \Phi^+$, $a_{e+1} \leq N_{e+1}<\langle \alpha, \lambda \rangle \leq \langle \tilde\alpha, \lambda \rangle = N_{e+2} < p-1-N_{e+1} \leq p-1-a_{e+1}$, where the second last inequality holds since $p>1+2N_{e+2}>1+N_{e+1}+N_{e+2}$. Thus the claim holds for $e+1$ in all cases. 
\end{proof}

Finally, we obtain the following theorem: 
\begin{thm} \label{thm:oneprime}
Let $G$ be a split reductive group with $\dim Z(G) \leq 1$. Let $e \geq 0$ and let $p$ be a prime number such that
\begin{enumerate}
    \item $p$ is larger than $n_{\alpha,\beta}$ for all $\alpha, \beta \in \Phi$, where $n_{\alpha, \beta}$ is defined in Proposition \ref{prop:chevalleybasis}.
    \item $p>1+2N_{e+1}$, where $N_{\bullet}$ is defined in Definition \ref{def:the integers N}.
    \item The index of irregularity $e_p$ is at most $e$, where $e_p$ is defined in Definition \ref{def:index}.
\end{enumerate}
Then there is a continuous representation $\rho \colon \galres{\Q}{\{p\}} \to G(\Zp)$ with open image. 
\end{thm}
\begin{proof}
If $G$ is semisimple, this follows immediately from Corollary \ref{cor:galrep one prime}. If $G$ is reductive with $\dim Z(G)=1$, the same argument as in the proof of Theorem \ref{thm:galrep} gives us the desired representation. 
\end{proof}

%\item $p$ is larger than $n_{\alpha,\beta}$ for all $\alpha, \beta \in \Phi$, where $n_{\alpha, \beta}$ is a positive integer depending only on the root system $\Phi$ (Proposition \ref{prop:chevalleybasis}).
    %\item $p>1+2N_{e+1}$, where $\{N_k\}$ is a sequence of integers depending only on the root system $\Phi$ (Definition \ref{def:the integers N}). 
    %\item The index of irregularity $e_p$ (Definition \ref{def:index}) is at most $e$. 

\end{document}